\documentclass[10pt]{amsart}
\usepackage{latexsym}
\usepackage{amssymb}
\usepackage{amsmath}
\usepackage{mathrsfs}
\usepackage[cp866]{inputenc}
\usepackage{bbm}
\usepackage{bbold}

\begin{document}

\title{
UPPER SEMILATTICES OF
FINITE-DIMENSIONAL GAUGES
}

\author{S.~S. Kutateladze}
\address[]{
Sobolev Institute of Mathematics\newline
\indent 4 Koptyug Avenue\newline
\indent Novosibirsk, 630090\newline
\indent Russia}
\email{
sskut@member.ams.org
}
\begin{abstract}
This is a brief overview of some applications of the
ideas of  abstract convexity to the
upper semilattices of gauges in finite dimensions.

\end{abstract}
\keywords
{Majorization, Minkowski functional, abstract convexity}
\dedicatory{In Memory of Alex Rubinov (1940--2006)}
\date{December 31, 2006}

\maketitle

\def\Aff{\operatorname{Aff}}
\def\ext{\operatorname{ext}}
\def\Fd{\operatorname{Fd}}
\def\spt{\operatorname{supp}}
\def\cop{\operatorname{cop}}
\def\ssection#1{\bigskip\begin{center}{\scshape#1}\end{center}\medskip}
 \def\endproc{\rm}
 \def\proclaim#1{{\bf #1}\sl}
 \def\Proclaim#1{\smallskip\par{\bf #1}\sl}
 \def\beginproof{{\sc Proof.}}
 \def\endproof{$\Box$}

\ssection{Introduction}

Duality in convexity is a simile  of reversal in positivity.
The ghosts of this similarity underlay the research on abstract
convexity we were engrossed in with Alex Rubinov in the early 1970s.
Our efforts led to the survey \cite{MD} and its expansion in the namesake
book~\cite{KutRub}. We always cherished a~hope to revisit this
area and shed light on a~few obscurities. However, the fate was
against us.

Inspecting the archive of our drafts of these years, I encountered several items
on the cones of Minkowski functionals or, equivalently, gauges.
The results on the Minkowski duality in finite dimensions are practically
unavailable in full form, whereas they rest on the technique that is still
uncommon and unpopular but definitely  profitable. The theorems on
gauges  appeared mostly in some mimeographed local sources that had disappeared
two decades ago. We hoped and planned to expatiate on these matters when time will come.

Alex Rubinov was my friend up to his terminal day. He shared  his
inspiration and impetus with me.  So does and will do his memory...

An abstract convex function is the upper envelope of a~family of simple
functions~\cite{MD}--\cite{Rub}. The cone of abstract convex elements is an upper semilattice.
We describe the bipolar of such a~semilattice through majorization generated
by its polar. Polyhedral approximation simplifies the generators of  the polar
in finite dimensions to discrete measures.
Decomposition reduces the matter to Jensen-type inequalities,
which opens a~possibility of linear programming and we are done.
These ideas  characterize our approach.

This article is organized as follows:  Section~1  is a~short discussion of majorization
and decomposition in the spaces of continuous functions.  Section~2
addresses the space of convex sets in finite dimensions and the influence  of
polyhedral approximation on the structure of  dual cones.  Section~3
illustrates the use of  linear programming for revealing continuous linear selections  over
convex figures.  Section~4   collects some dual representations
for the members of upper semilattices of gauges. In Section~5 we deal with
some upper lattices of gauges that are closed under intersection.

\section{Majorization and Decomposition}

It was long ago in 1954 that  Reshetnyak suggested
in his unpublished  thesis \cite{Resh} to compare (positive) measures on the
Euclidean unit sphere $S_{N-1}$ as follows:

\subsection{}
 A measure $\mu$ {\it (linearly) majorizes} or {\it dominates}
a~measure $\nu$ provided that to each decomposition of
$S_{N-1}$ into finitely many disjoint Borel sets $U_1,\dots,U_m$
there are measures $\mu_1,\dots,\mu_m$ with sum $\mu$
such that every difference $\mu_k - \nu|_{U_k}$ annihilates
all restrictions to $S_{N-1}$ of linear functionals over
$\mathbbm R^N$. In symbols, we write $\mu\gg{}_{\mathbbm R^N} \nu$.

Reshetnyak  proved that
$$
\int\limits_{S_{N-1}} p d\mu \ge  \int\limits_{S_{N-1}} p d\nu
$$
for every  sublinear
functional  $p$
on  $\mathbbm R^N$   if   $\mu\gg{}_{\mathbbm R^N} \nu$.
This gave an important trick for generating positive linear functionals
over various classes of convex  surfaces and functions.

\subsection{}
A~similar idea was suggested by
Loomis \cite{Loomis} in 1962 within Choquet theory.
A~measure $\mu$ {\it affinely majorizes}
a~measure $\nu$, both given
on a~compact convex subset $Q$ of a~locally convex space $X$,
provided that
to each decomposition of
$\nu$ into finitely many summands
$\nu_1,\dots,\nu_m$  there are measures $\mu_1,\dots,\mu_m$
with $\mu$ such that every difference
$\mu_k - \nu_k$ annihilates all restrictions
to   $Q$  of the affine   functions over $X$.
In symbols, $\mu\gg{}_{\Aff(Q)} \nu$.
Many applications of affine majorization are set forth in~\cite{MarOl}.

Cartier, Fell, and Meyer  proved in~\cite{CFM} that
$$
\int\limits_{Q} f d\mu \ge  \int\limits_{Q} f d\nu
$$
for every continuous convex function  $f$
on  $Q$   if and only if   $\mu\gg{}_{\Aff(Q)} \nu$.

An analogous necessity part for linear majorization was published
in~\cite{KutDAN}. In applications we use a~more
detailed version of majorization~\cite{KutCho}:

\subsection{}
\proclaim{Decomposition Theorem.}
Assume that $H_1,\dots,H_n$ are  cones in a~vector lattice~$X$.
Assume further that $f$ and $g$ are positive linear functionals on~$X$.
The inequality
$$
f(h_1\vee\dots\vee h_n)\ge g(h_1\vee\dots\vee h_n)
$$
holds for all
$h_k\in H_k$ $(k:=1,\dots,n)$
if and only if to each decomposition
of~ $g$ into a~sum of~$n$ positive terms
$g=g_1+\dots+g_N$
there is a~decomposition of ~$f$ into a~sum of~$n$
positive terms $f=f_1+\dots+f_n$
such that
$$
f_k(h_k)\ge g_k(h_k)\quad
(h_k\in H_k;\ k:=1,\dots,n).
$$
\endproc

\section{The Space of Convex Figures}

We will proceed in the Euclidean space $\mathbbm R^N$.

\subsection{}
A~{\it convex figure\/} is a~compact convex set. A~{\it convex body\/}
is a~solid convex figure.
The {\it Minkowski  duality\/} identifies
a~convex figure $S$ in
$\mathbbm R^N$ with its {\it support function\/}
$S(z):=\sup\{(x,z)\mid  x\in S\}$ for $z\in \mathbbm R^N$.
Considering the members of $\mathbbm R^N$ as singletons, we assume that
$\mathbbm R^N$ lies in the set $\mathscr V_N$
of all compact convex subsets
of $\mathbbm R^N$.

\subsection{}
The Minkowski duality makes $\mathscr V_N$ into a~cone
in the space $C(S_{N-1})$  of continuous functions on the Euclidean unit sphere
$S_{N-1}$, the boundary of the unit ball $\mathfrak z_N$.
This yields  is  the so-called {\it  Minkowski structure\/} on $\mathscr V_N$.
Addition of the support functions
of convex figures amounts to taking their algebraic sum, also called the
{\it Minkowski addition}. It is worth observing that the
{\it linear span}
$[\mathscr V_N]$ of~$\mathscr V_N$ is dense in $C(S_{N-1})$, bears
a~natural structure of a~vector lattice
and is usually referred to as the {\it space of convex sets}.
The study of this space stems from the pioneering breakthrough of
Alexandrov~\cite{AD} in 1937 and the further insights of  Radstr\"{o}m~\cite{Rad}
and H\"{o}rmander~\cite{Her}.

\subsection{}
A~{\it gauge} $p$ is a~positive sublinear functional on a~real vector space~
$X$ viewed as the Minkowski
functional of the  conic segment  $S_p:=\{p\le1\}:=\{x\in X\mid p(x)\le 1\}$.
The latter is also referred to as a~{\it gauge} or {\it caliber}.
A~gauge $p$ is a~{\it norm\/} provided that its {\it ball\/} $S_p$
is symmetric and absorbing.
Recall that the {\it subdifferential\/}
or  {\it support set\/} $\partial p$ of $p$
is the {\it dual ball\/} or {\it polar\/} of $S_p$. The polar of a~ball
$S$ is denoted by $S^\circ$ and the dual norm of $\|\cdot\|_S$
is $\|\cdot\|_{S^{\circ}}$.
The ``donkey bridge'' of functional analysis consists in the
{\it duality rules}:
$$
\|\cdot\|_S =S^{\circ}(\cdot),\quad  \|\cdot\|_{S^{\circ}}=S(\cdot).
$$
We will restrict exposition  to the norms  and
balls of~$\mathbbm R^N$ by way of tradition.

\subsection{}
\proclaim{Approximation Lemma.} If $H$ is a subcone of $\mathscr V_N$ then
the signed measures with finite support are sequentially weakly* closed
in the dual cone $H^{*}$.
\endproc

{\sc Proof}.
Let $\mu\in H^{*}$.
The mappings
$$
z\mapsto \mu_{+}(z);\quad
z\mapsto \mu_{-}(z),
$$
with $z\in\mathbbm{R}^{N}$, are linear functionals on~$\mathbbm{R}^{N}$.
Therefore, there are $u,v\in\mathbbm{R}^{N}$ such that
$\mu_{+}(z)=(u,z)$ and $\mu_{-}(z)=(v,z)$.
Put
$$
\gathered
\overline{\mu}_{1};=\mu_{+}+\text{mes}+|u|\varepsilon_{-u/|u|};
\\
\overline{\mu}_{2}:=\mu_{-}+\text{mes}+|v|\varepsilon_{-v/|v|};
\\
\mu_{1}:=\overline{\mu}_{1}+|v|\varepsilon_{-v/|v|};
\quad
\mu_{2}:=\overline{\mu}_{2}+|u|\varepsilon_{-u/|u|}.
\endgathered
$$
As usual, $\varepsilon_z$  is the Dirac measure at $z\in \mathbbm R^N$, while
$|\cdot|$ is the Euclidean norm on $\mathbbm R^N$, and  mes is the Lebesgue measure on~$S_{N-1}$:
i.~e. the surface area function of the
Euclidean ball
$\mathfrak{z}_{N}:=\{x\in\mathbbm{R}^{N}:|x|\le1\}$.
Note that  $\mu=\mu_{1}-\mu_{2}$. Moreover,
the measures $\overline{\mu}_{1}$ and $\overline{\mu}_{2}$
are nondegenerate and translation-invariant.
Indeed, check that so is $\overline{\mu}_{1}$.
This signed measure is clearly positive and not supported by any great
hypersphere. We  are left with validating translation-invariance.
If $k:=1,\dots,N$ then
$$
\int\limits_{S_{N-1}}e_{j}d\overline{\mu}_{1}=
\int\limits_{S_{N-1}}e_{j}d\mu_{+}+
\int\limits_{S_{N-1}}e_{j}d\mu(\mathfrak z_{N})-
(u,e_{k})=(u,e_{k})-(u,e_{k})=0.
$$
Consider  a convex figure $\mathfrak{x}$ whose
surface area function  $\mu(\mathfrak{x})$ equals $\overline{\mu}_{1}$.
The existence of this figure is guaranteed by the celebrated
Alexandrov Theorem~\cite[p.108]{AD}.

Let $(\mathfrak{x}_{m})$ be a sequence of polyhedra including
$\mathfrak{x}$ and converging to $\mathfrak{x}$ in the
Hausdorff metric on $[\mathscr V_N]$which is induced by the
Chebyshev norm on $C(S_{N-1})$.
Then the measures $\overline{\mu}_{m}^{1}=\mu(\mathfrak{x}_{m})$ converge
weakly* to $\overline{\mu}_{1}$ and
 $\overline{\mu}_{m}^{1}\gg{}_{\mathbbm R^n}\overline{\mu}_{1}$.
Indeed, given a convex figure $\mathfrak{z}$, we have
$$
\gathered
\int\limits_{S_{N-1}}\mathfrak{z}
d\overline{\mu}_{m}^{1}=
\int\limits_{S_{N-1}}\mathfrak{z}
d\mu(\mathfrak{x}_{m})=nV(\mathfrak{z},\mathfrak{x}_{m},\dots,\mathfrak{x}_{m})
\\
\ge nV(\mathfrak{z},\mathfrak{x},\dots,\mathfrak{x})=
\int\limits_{S_{N-1}}\mathfrak{z}
d\mu(\mathfrak{x})=\int\limits_{S_{N-1}}\mathfrak{z}
d\overline{\mu}_{1}
\endgathered
$$
by the inclusion monotonicity of the mixed volume $V(\cdot,\dots,\cdot)$
in every argument..
By analogy, there is a sequence $(\overline{\mu}_{m}^{2})$,
converging weakly* to $\overline{\mu}_{2}$  and such that
$\overline{\mu}_{m}^{2}\gg{}_{\mathbbm R^n}\overline{\mu}_{2}$.
Putting
$$
\gathered
\mu_{m}^{1}:=\overline{\mu}_{m}^{1}+|v|
\varepsilon_{-v/|v|};
\\
\mu_{m}^{2}:=\overline{\mu}_{m}^{2}+|u|
\varepsilon_{-u/|u|},
\endgathered
$$
we see that $\mu_{m}^{1}-\mu_{m}^{2}$ converges weakly* to~$\mu$.
The proof is complete.

\section{Labels and Decompositions}

The Approximation Lemma allows us to reduce consideration
to signed measures with finite support.  These measures
decompose easily. We will exhibit a typical application.

\subsection{}
A family $(\mu_1,\dots,\mu_n)$
of regular Borel measures on the sphere  $S_{N-1}$
is a~{\it labeling\/} on $\mathbbm R^N$  provided that
$(\mu_1(\mathfrak{x})$, $\dots$, $\mu_n(\mathfrak{x}))\in \mathfrak{x}$
for all $\mathfrak{x}\in\mathscr{V}_N$.
The vector   $(\mu_1(\mathfrak{x}), \dots,\mu_n(\mathfrak{x}))$
is a~{\it label\/} of~$\mathfrak{x}$.

\subsection{}
\proclaim{Proposition.}
A~family  $(\mu_1,\dots,\mu_n)$
is a labeling on $\mathbbm R^N$ if and only if
$$
\varepsilon_x-\sum\limits_{k=1}^n x_k\mu_k\in\mathscr{V}_{N}^*.
$$
for all $x\in S_{N-1}$.
\endproc

{\sc Proof.}
The Minkowski duality is an isomorphism of the relevant structures. Hence, the
definition of labeling can be rephrased as follows:
$$
\sum\limits_{k=1}^n x_k\mu_k(\mathfrak{x})\le \mathfrak{x}(x)\quad
(x\in \mathbbm{R}^n,\ \mathfrak{x}\in\mathscr{V}_N).
$$

\subsection{}
Using linear majorization for describing  $\mathscr{V}_{N}^*$,
we can write down some criteria for labeling in terms of decompositions.
For simplicity, we will argue in the planar case.

Consider the conditions:
$$
\leqno{(++)}\qquad
\varepsilon_{(\vartriangle_1,\vartriangle_2)}+\vartriangle_1\mu_1^-
+\vartriangle_2\mu_2^-
\underset{\mathbbm{R}^2}{\gg}\ \vartriangle_1\mu_1^+
+\vartriangle_2\mu_2^+;
$$
$$
\leqno{(+-)}\qquad
\varepsilon_{(\vartriangle_1,-\vartriangle_2)}+\vartriangle_1\mu_1^-
+\vartriangle_2\mu_2^+
\underset{\mathbbm{R}^2}{\gg}\ \vartriangle_1\mu_1^+
+\vartriangle_2\mu_2^-;
$$
$$
\leqno{(-+)}\qquad
\varepsilon_{(-\vartriangle_1,\vartriangle_2)}+\vartriangle_1\mu_1^+
+\vartriangle_2\mu_2^-
\underset{\mathbbm{R}^2}{\gg}\ \vartriangle_1\mu_1^-
+\vartriangle_2\mu_2^+;
$$
$$
\leqno{(--)}\qquad
\varepsilon_{(-\vartriangle_1,-\vartriangle_2)}+\vartriangle_1\mu_1^+
+\vartriangle_2\mu_2^+
\underset{\mathbbm{R}^2}{\gg}\ \vartriangle_1\mu_1^-
+\vartriangle_2\mu_2^-;
$$
with $(\vartriangle_1,\vartriangle_2)\in S_1\cap\mathbbm{R}^2_+$.
Clearly, the requirement of  4.1 amounts to the four conditions
simultaneously.
By way of example, we will elaborate the relevant criterion only
in the case of~$(+-)$.

\subsection{}
\proclaim{Proposition.}
For $(+-)$ to hold it is necessary and sufficient that to
all $(\vartriangle_{1},\vartriangle_{2})$ in $S_{1}\cap\mathbbm{R}_{+}^{2}$ and all decompositions
$\{(\mu_{1}^{+})_1,\dots,(\mu_{1}^{+})_m\}$
of ~$\mu_{1}^{+}$
and al decompositions
$\{(\mu_{2}^{-})_1,\dots,(\mu_{2}^{-})_m\}$
 of $\mu_{2}^{-}$
there exist a~decomposition  $\{(\mu_{1}^{-})_1,\dots,(\mu_{1}^{-})_m\}$
of~$\mu_{1}^{-}$, a~decomposition  $\{(\mu_{2}^{+})_1,\dots,(\mu_{2}^{+})_m\}$
of~$\mu_{2}^{+}$, and reals $\alpha_{1},\dots,\alpha_{m}$
that make compatible the simultaneous inequalities:
$$
\gathered
\alpha_{1}\geq0;\dots;\alpha_{m}\geq0;
\alpha_{1}+ \ldots +\alpha_{m}=1;
\\
\vartriangle_{1}(x_{(\mu_{1}^{-})_{k}}-x_{(\mu_{1}^{+})_k}
+\alpha_{k}e_{1})
\\
=\vartriangle_{2}(x_{(\mu_{2}^{-})_{k}}-x_{(\mu_{2}^{+})_k}
+\alpha_{k}e_{2})\quad (k:=1,\dots,m),
\endgathered
$$
where $x_{\mu}$ is the representing point of~ $\mu$; i.~e.,
$\mu(u)=(u,x_{\mu})$ for all $u\in\mathbbm{R}^{2}$.
\endproc

{\sc Proof.}
$\Longleftarrow$:
Let $(\vartriangle_{1},\vartriangle_{2})\in S_{1}\cap\mathbbm{R}_{+}^{2}$
and let $\{\nu_{1},\dots,\nu_{m}\}$
be an~arbitrary decomposition of~
$\vartriangle_{1}\mu_{1}^{+}+\vartriangle_{2}\mu_{2}^{-}$.
By the Riesz Decomposition Lemma there are a~decomposition
$\{(\mu_{1}^{+})_{1},\dots,(\mu_{1}^{+})_{m}\}$
of~$\mu_{1}^{+}$ and a~decomposition
$\{(\mu_{2}^{-})_{1},\dots,(\mu_{2}^{-})_{m}\}$
of~$\mu_{2}^{-}$ such that
$\vartriangle_{1}(\mu_{1}^{+})_k+\vartriangle_{2}(\mu_{2}^{-})_k=\nu_k$.
Find some parameters satisfying the simultaneous inequalities and put
$$
\mu_{k}:=\vartriangle_{1}(\mu_{1}^{-})_k+
\vartriangle_{2}(\mu_{2}^{+})_k+\alpha_{k}
\varepsilon_{(\vartriangle_{1},-\vartriangle_{2})}.
$$
Clearly, $\mu_{k}\geq0$ and, moreover,
$$
\sum\limits_{k=1}^{m}\mu_{k}=\vartriangle_{1}\mu_{1}^{-}+
\vartriangle_{2}\mu_{2}^{+}+\varepsilon_{(\vartriangle_{1},-\vartriangle_{2})}.
$$
Furthermore,
$$
\gathered
x_{\mu_{k}}-x_{\nu_{k}}=\vartriangle_{1} x_{(\mu_{1}^{-})_k}
+\vartriangle_{2}x_{(\mu_{2}^{+})_k}+\alpha_{k}\vartriangle_{1}e_{1}
\\
-\alpha_{k}\vartriangle_{2}e_{2}-
\vartriangle_{1}x_{(\mu_{1}^{+})_k}-\vartriangle_{2}x_{(\mu_{2}^{-})_k}=0,
\endgathered
$$
and so $\mu_{k}-\nu_{k}$ belongs to the polar of~$\mathbbm{R}^{2}$ in~$C(S_{1})$.

$\Longrightarrow$:
Assume $(+-)$ valid.

Given decompositions
$\{(\mu_{1}^{+})_{1},\dots,(\mu_{1}^{+})_{m}\}$ and
$\{(\mu_{2}^{-})_{1},\dots,(\mu_{2}^{-})_{m}\}$
there is a~decomposition $\{\nu_{1},\dots,\nu_{2m}\}$
of~$\varepsilon_{(\vartriangle_{1},-\vartriangle_{2})}
+\vartriangle_{1}\mu_{1}^{-}+\vartriangle_{2}\mu_{2}^{+}$
such that
$$
x_{\nu_{k}}=x_{(\mu_{1}^{+})_k};\quad
x_{\nu_{m+k}}=x_{(\mu_{2}^{-})_k}\quad
(k:=1,\dots,m).
$$

We are left with appealing to the Riesz Decomposition Lemma
and representing the decomposition
$\{\nu_{1},\dots,\nu_{2m}\}$
through the corresponding decompositions of~$\varepsilon_{(\vartriangle_{1},-\vartriangle_{2})}$,
$\vartriangle_{1}\mu_{1}^{-}$, and
$\vartriangle_{2}\mu_{2}^{+}$.
The proof is complete.

\subsection{}
If it is possible to chose decompositions in 3.4  independently
of~$(\vartriangle_{1},\vartriangle_{2})$,
then we  come to a~sufficient condition for labeling.
Let us illustrate this by exhibiting an example of one of the simplest labelings.

We will seek a labeling of the form
$$
\mu_{1}:=|\mu^{+}|\varepsilon_{\mu^{+}/|\mu^{+}|}
-|\mu^{-}|\varepsilon_{\mu^{-}/|\mu^{-}|};
$$
$$
\mu_{2}:=|\nu^{+}|\varepsilon_{\nu^{+}/|\nu^{+}|}
-|\nu^{-}|\varepsilon_{\nu^{-}/|\nu^{-}|},
$$
with $\mu^{+}$, $\mu^{-}$, $\nu^{+}$, and  $\nu^{-}$ some points on the plane.
The sufficient condition we have just suggested paraphrases as follows:
$$
\gathered
\alpha_{k},\beta_{k},a_{k},b_{k},\gamma_{k},c_{k}\geq0;
\\
\alpha_{k}+a_{k}=1;\quad
\beta_{k}+b_{k}=1;\quad
\gamma_{k}+c_{k}=1\quad (k:=1,\dots,4);
\\
\mu^{+}=\alpha_{1}\mu^{-}+\gamma_{1}e_{1};\quad
\beta_{1}\nu^{-}+\gamma_{1}e_{2}=0;
\\
\nu^{+}=b_{1}\nu^{-}+c_{1}e_{2};\quad
a_{1}\mu^{-}+c_{1}e_{1}=0;
\\
\mu^{-}=\alpha_{2}\mu^{+}-\gamma_{2}e_{1};\quad
\beta_{2}\nu^{-}+\gamma_{2}e_{2}=0;
\\
\nu^{+}=b_{2}\nu^{-}+c_{2}e_{2};\quad
a_{2}\mu^{+}-c_{2}e_{1}=0;
\\
\mu^{-}=\alpha_{3}\mu^{+}-\gamma_{3}e_{1};\quad
\beta_{3}\nu^{+}-\gamma_{3}e_{2}=0;
\\
\nu^{-}=b_{3}\nu^{+}-c_{3}e_{2};\quad
a_{3}\mu^{+}-c_{3}e_{1}=0;
\\
\mu^{+}=\alpha_{4}\mu_{-}+\gamma_{4}e_{1};\quad
\beta_{4}\nu^{+}-\gamma_{4}e_{2}=0;
\\
\nu^{-}=b_{4}\nu^{+}-c_{4}e_{2};\quad
a_{4}\mu^{-}+\gamma_{4}e_{1}=0.
\endgathered
$$
The solution of the last system is given by the parameters:
$$
\gathered
\alpha_{k}=b_{k}=0;\quad
\beta_{k}=a_{k}=1;
\\
\gamma_{k}=c_{k}=\frac{1}{2}\quad (k:=1,\dots,4).
\endgathered
$$
Moreover,
$$
\mu^{+}=\frac{1}{2}e_{1};\quad \nu^{+}=\frac{1}{2}e_{2};
\quad
\mu^{-}=-\frac{1}{2}e_{1};\quad \nu^{-}=-\frac{1}{2}e_{2}.
$$
Therefore, the simplest labeling of~$\mathfrak{x}$
is the point
$\frac{1}{2}(\mathfrak{x}(e_{1})-\mathfrak{x}(-e_{1}),
\mathfrak{x}(e_{2})-\mathfrak{x}(-e_{2}))$.
It is worth emphasizing that the validation of the above conditions belongs
to linear programming which enables us to seek for arbitrary labelings by
signed measures with finite support.

\section{The Case of Joining Gauges}

We now apply the above ideas to studying the classes of
$N$-dimensional convex surfaces which comprise upper semilattices
in~$\mathscr{V}_N$. To simplify notation we will discuss only
balls, denoting the set of  balls in
$\mathscr{V}_N$ by $\mathscr VS_N$. It is convenient formally
to add the apex to $\mathscr VS_N$.
If $S\in\mathscr{V}_N S$ differs from the origin then we
use the symbol $\|\cdot\|_S$ not only for the gauge of~$S$
but also for the {\it operator norm\/} corresponding to $S$
in the endomorphism space
$\mathscr{L}(\mathbbm{R}^{N})$ of~$\mathbbm{R}^{N}$.
In other words,
$$
\|x\|_{S}:=\inf\{\alpha>0\mid x/\alpha\in S\}\quad (x\in\mathbbm{R}^{N});
$$
$$
\|A\|_{S}:=\sup\{\|Ax\|_{S}\mid x\in S\}\quad
(A\in\mathscr{L}(\mathbbm{R}^{N})).
$$
Recall that
$$
S^\circ=\{x\in\mathbbm{R}^{N}\mid |(x,y)|\le1\ (y\in S)\},
$$
where $(\cdot,\cdot)$ is the standard inner product of~$\mathbbm R^N$.

Observe that $\mathscr{V}_{N}S$
is a lattice and simultaneously a~cone. However,
$\mathscr{V}_{N}S$ is not closed in~$\mathscr{V}_{N}$.
This circumstance notwithstanding,  given a~family
$(S_{\xi})_{\xi\in\Xi}$ in~$\mathscr{V}_{N}S$,
sometimes
we may soundly speak of the {\it upper hull\/}
$\pi^\uparrow(\Xi)$, {\it lower hull\/} $\pi_\downarrow(\Xi)$, and {\it hull} $\pi(\Xi)$
of this family, implying the least closed cones that
lie  in~$\mathscr{V}_{N}S$, include $S_{\xi}$ for all $\xi\in\Xi$,
and are closed under the join, the meet, and both operations in the
lattice of convex figures $\mathscr{V}_{N}$.
An example  is provided by any instance of {\it nondegenerate family}. The latter is by definition
any family of nonzero sets $(S_{\xi})_{\xi\in\Xi}$ such that,
$$
\sup\limits_{\xi\in\Xi}\|A\|_{S_{\xi}}<+\infty\quad
(A\in\mathscr{L}(\mathbbm{R}^{N})).
$$
Indeed, put
$$
\mathscr{A}(\Xi):=\{A\in\mathscr{L}(\mathbbm{R}^{N})\mid AS_{\xi}\subset S_{\xi}\ (\xi\in\Xi)\},
$$
and let ${\rm M}(\Xi)$ be the set of the symmetric elements of~$\mathscr{V}_N$ such that $AS\subset S$
for all $A\in\mathscr{A}(\Xi)$.
Since  $(S_{\xi})_{\xi\in\Xi}$ is nondegenerate, all
members of ${\rm M}(\Xi)$ but the zero singleton are absorbing.
Moreover, ${\rm M}(\Xi)$ is clearly a closed sublattice of~$\mathscr{V}_N$.

We will need the helpful property of a~nondegenerate family:
If $y\in\mathbbm{R}^{N}$ differs from the zero~of~$\mathbbm{R}^{N}$ then
$$
S_{y}:=\bigwedge\limits_{\xi\in\Xi}\frac{S_{\xi}}{S_{\xi}(y)}
$$
is absorbing. Indeed, given  $z\in\mathbbm{R}^{N}$ we infer that
$$
\sup\limits_{\xi\in\Xi}S_{\xi}(y)S_{\xi}^{\circ}(z)=
\sup\limits_{\xi\in\Xi}\|y\|_{S_{\xi}^{\circ}}
\|z\|_{S_{\xi}}=
\sup\limits_{\xi\in\Xi}\|y\otimes z\|_{S_{\xi}}<+\infty,
$$
where $y\otimes z:x\mapsto (y,x)z$ for all $x\in \mathbbm R^N$.
Hence, the polar of~$S_{y}$ is compact, which implies that
$S_{y}$ is absorbing. Without further specification,
we will address only nondegenerate families of balls in the sequel.

\subsection{}
\proclaim{Theorem.}
A~gauge $S$ belongs to~$\pi^\uparrow(\Xi)$
if and only if
$$
\frac{S}{\sum\limits_{k=1}^{n}\|x_{k}\|_{S^{\circ}}}\le
\bigvee\limits_{\xi\in\Xi}
\frac{S_{\xi}}{\sum\limits_{k=1}^{n}\|x_{k}\|_{S^{\circ}_{\xi}}}
$$
for any collection of the vectors $x_{1},\dots,x_{p}\in\mathbbm{R}^{N}$
that are not all zero simultaneously.
\endproc

{\sc Proof.}
It is obvious that $\pi^\uparrow(\Xi)$ is the closure of the upper semilattice
of all  $H$-convex functions
with $H$ the conic hull of the family $(S_{\xi})_{\xi\in\Xi}$.
The polar of~$\pi^\uparrow(\Xi)$ may be approximated
with finitely supported signed measures by the Approximation Lemma.
Using the Bipolar Theorem, we see that
$S\in\pi^\uparrow(\Xi)$
if and only if
$\sum\nolimits_{k=1}^{n}S(x_{k})\ge S(y)$
whenever $y,x_{1},\dots,x_{n} \in \mathbbm{R}^{N}$
satisfy $\sum\nolimits_{k=1}^{n}S_{\xi}(x_{k})\ge S_{\xi}(y)$
for all $\xi\in \Xi$.
By  duality, $S\in\pi^\uparrow(\Xi)$
if and only if
$$
\bigwedge\limits_{\xi\in\Xi}\sum\limits_{k=1}^{n}
\|x_{k}\|_{S^{\circ}_{\xi}}S^{\circ}_{\xi}
\subset\sum\limits_{k=1}^{n}\|x_{k}\|_{S^{\circ}}S^{\circ}.
$$
Taking polars, we complete the proof of the theorem.

\subsection{}
\proclaim{Corollary.}
A nonzero gauge  $S$ belongs to $\pi^\uparrow(\Xi)$
if and only if
$$
S=\bigwedge\limits_{(x_{1},\dots,x_{n})}
\sum\limits_{k=1}^{n}S(x_{k})
\bigvee\limits_{\xi\in\Xi}
\frac{S_{\xi}}{\sum\limits_{k=1}^{n}S_{\xi}(x_{k})},
\leqno{(4.2.1)}
$$
where the intersection ranges over all nonzero tuples
$(x_{1},\dots,x_{n})\in (\mathbbm{R})^{N}$.
\endproc

{\sc Proof.}
Clearly, (4.2.1) guarantees the inclusion of 4.1
and so $S\in\pi^\uparrow(\Xi)$.
The last containment in turn implies
the simple representation:
$$
S=\bigwedge\limits_{x\ne 0}S(x)\bigvee\limits_{\xi\in\Xi}
\frac{S_{\xi}}{S_{\xi}(x)}.
\leqno{(4.2.2)}
$$
Indeed, denote by~$\widetilde{S}$
the right-hand side of~(4.2.2).
By 4.1, $S\le\widetilde{S}$.
If $z\in\mathbbm{R}^{n}$ then
$$
\gathered
\widetilde{S}(z)=\biggl(\,\bigwedge\limits_{x\ne0}
S(x)\bigvee\limits_{\xi\in\Xi}\frac{S_{\xi}}{S_{\xi}(x)}
\biggr)(z)
\\
\le S(z)\biggl(\,\bigvee\limits_{\xi\in\Xi}
\frac{S_{\xi}}{S_{\xi}(z)}\biggr)(z)=
S(z)\bigvee\limits_{\xi\in\Xi}
\frac{S_{\xi}(z)}{S_{\xi}(z)}=S(z).
\endgathered
$$
By the Minkowski duality $\widetilde{S}\le S$.
Denote by $\overset{\approx}{S}$
the right-hand side of~(4.2.1).
Since $S\le\overset{\approx}{S}\le\widetilde{S}\le S$;
therefore, $S=\overset{\approx}{S}$ and we are done.

\subsection{}
From 4.2 it follows that if each closed subset
of~$\mathscr{V}_{n}S$ is a cone provided that it contains the convex hull
and intersection
of any pair of its elements as well as the
dilation $\alpha\mathfrak{x}$, with $\alpha\ge0$,
of its every member $\mathfrak{x}$.

\subsection{}
The proof of Theorem 4.1 shows that
a positively homogeneous continuous function~$f$ on~$\mathbbm R^N$
 is the support function of a~member of~$\pi^\uparrow(\Xi)$ if and only if
$\sum\limits_{k=1}^{n}f(x_k)\ge f(y)$
provided that
$\sum\limits_{k=1}^{n}S_{\xi}(x_k)\ge S_{\xi}(y)$
for all $\xi\in\Xi$.
Observe that we may restrict the range of the index to~$n=1$
only on condition that the balls $S_{\xi}$
are dilations of one another.
Indeed, in this event the polar $\pi^\uparrow(\Xi)$
is the weakly* closed conic hull of two-points relations
and so the functions of the form
$x\mapsto\alpha S_{\xi_1}(x)\wedge\beta S_{\xi_2}(x)$
turn out sublinear for positive $\alpha$ and~$\beta$.

\section{The Case of Meeting Gauges}
We now address some properties of gauges which are tied with
intersection. This operation involves some peculiarities since the intersection
of balls differs in general from the pointwise infimum of their
support functions. However, the idea of decomposition
applies partially to this case.

\subsection{}
\proclaim{Theorem.}
Let $H$ be a~cone in~$\mathscr{V}_{N}S$
and $H=\pi_\downarrow (H)$.
Assume given  a nonzero vector $y$ in~$\mathbbm{R}^{N}$
such that
$$
S_{y}:=
\bigvee\limits_{S\in H;S\ne\{0\}}
\frac{S}{S(y)}
$$
is absorbing. Take $x_{1},\dots,x_{n}$ in $\mathbbm{R}^{N}$.
The inequality
$$
\sum\limits_{k=1}^{n}S(x_{k})\ge S(y)
$$
holds for every gauge $S\in H$
if and only if there are vectors
$z_{1},\dots,z_{n}$ in $\mathbbm{R}^{N}$ such that
$\sum\limits_{k=1}^{n}z_{k}=y$ and,
moreover, $S(x_{k})\ge S(z_{k})$ for all~$S\in H$.
\endproc

{\sc Proof.}
$\Longleftarrow$:
Since $S$ is a gauge, the support function
of~$S$ is a~sublinear functional and
$$
\sum\limits_{k=1}^{n}S(x_{k})\ge
\sum\limits_{k=1}^{n}S(z_{k})\ge
S\biggl(\,\sum\limits_{k=1}^{n}z_{k}\biggr)=S(y).
$$

$\Longrightarrow$:
For simplicity we restrict exposition to the case when
$S_{y}$ is absorbing for every nonzero $y\in\mathbbm{R}^{N}$.
Put
$$
K:=\sup\limits_{x\in S_{y}^{\circ}}|x|.
$$
By hypotheses, $K<+\infty$.
We further put
$$
\gathered
U:=\Bigl\{(\nu_{1},\nu_{2})\in C'(S_{N-1})\times
C'(S_{N-1})\mid\nu_{1}\ge0,\nu_{2}\ge0;
\\
\|\nu_{1}\|\vee\|\nu_{2}\|\le K;\quad
\int\limits_{S_{N-1}}(l,\cdot)
d(\nu_{1}+\nu_{2})=(l,y)\quad (l\in\mathbbm{R}^{N})\Bigr\};
\\
\widetilde{U}:=U+H^*\times H^*;
\\
\mu_{1}:=|x_{1}|\varepsilon_{x_{1}/|x_{1}|};\quad
\mu_{2}:=\sum\limits_{k=2}^{n}|x_{k}|\varepsilon_{x_{k}/|x_{k}|}.
\endgathered
$$
As usual, we agree that  the symbol $|0|\varepsilon_{0/|0|}0$ stands for the zero
vector.

Assume that the pair $(\mu_{1},\mu_{2})$
does not belong to~$\widetilde{U}$. Since $U$
is a weakly* compact convex set; therefore,
$\widetilde{U}$ is weakly* closed and convex.
By the Separation Theorem there are nonzero functions
$S'_{1}$ and $S'_{2}$ in~$H$
such that
$$
\mu_{1}(S'_{1})+\mu_{2}(S'_{2})<\nu_{1}(S'_{1})+\nu_{2}(S'_{2})
\leqno{(5.1.1)}
$$
for all $(\nu_{1},\nu_{2})\in U$.
Put
$$
S_{1}:=\frac{S'_{1}}{S'_{1}\wedge S'_{2}(y)};\quad
S_{2}:=\frac{S'_{2}}{S'_{1}\wedge S'_{2}(y)}.
$$
Note that $S_{1},S_{2}\in H$.
Consequently, the meet $S_{1}\wedge S_{2}$
belongs to~$H$. Moreover,
$$
\|y\|_{S_{1}^{\circ}\vee S_{2}^{\circ}}=
(S_{1}^{\circ}\vee S_{2}^{\circ})^{\circ}(y)=
S_{1}\wedge S_{2}(y)=\frac{S'_{1}\wedge S'_{2}}{S'_{1}\wedge S'_{2}(y)}
(y)=1.
$$
Since
$S_{1}\wedge S_{2}\supset S_{y}$; therefore,
 $S_{1}^{\circ}\vee S_{2}^{\circ}\subset S_{y}^{\circ}$.
In particular,
$$
\sup\limits_{x\in S_{1}^{\circ}\vee S_{2}^{\circ}}|x|\le K
\leqno{(5.1.2)}
$$
Let $V$ be a~face of~$S_{1}^{\circ}\vee S_{2}^{\circ}$
that contains~$y$; i.~e.,  the intersection of~$S_{1}^{\circ}\vee S_{2}^{\circ}$
with  some supporting hyperplane to~$S_{1}^{\circ}\vee S_{2}^{\circ}$
at~$y$.
Denote by  $\ext(V)$ the set of extreme points of~$V$.
By the Choquet Theorem there is a~probability measure $\overline{\nu}$
with support $\ext(V)$ and barycenter~$y$.
Put $V_{1}:=\ext(V)\cap S_{1}^{\circ}$ and
$V_{2}:=\ext(V)\setminus V_{1}$.
The set $V_{2}$ lies in~$S_{2}^{\circ}$.
Let $\overline{\nu}_{1}:=\overline{\nu}\vert_{V_{1}}$
and $\overline{\nu}_{2}:=\overline{\nu}\vert_{V_{2}}$.
Then $\overline{\nu}=\overline{\nu}_{1}+\overline{\nu}_{2}$.

We will treat a continuous function $f$ on~$S_{N-1}$
as the restriction to $S_{N-1}$ of the unique
positively homogeneous namesake function on~$\mathbbm{R}^{N}$ and put
$$
\nu_{1}:f\mapsto\int\limits_{V_{1}}fd\overline{\nu}_{1};
$$
$$
\nu_{2}:f\mapsto\int\limits_{V_{2}}fd\overline{\nu}_{2}\quad
(f\in C(S_{N-1}));
$$
$$
\nu:=\nu_{1}+\nu_{2}.
$$
Using (5.1.2) and the estimate $\overline{\nu}_{1}(\mathbbm{1})\le
\overline{\nu}(\mathbbm{1})=1$, with $\mathbbm 1$ the identically one function; we see that
$$
\|\nu_{1}\|=\nu_{1}(\mathbbm{1})=\int\limits_{V_{1}}
|\cdot|d\overline{\nu}_{1}\le
\sup\limits_{x\in S_{1}^{\circ}\vee S_{2}^{\circ}}|x|<K.
$$
\noindent
By analogy $\|\nu_{2}\|\le K$.
Moreover,
$$
\nu(l)=\int\limits_{V_{1}}(l,\cdot)d\overline{\nu}_{1}+
\int\limits_{V_{2}}(l,\cdot)d\overline{\nu}_{2}=
\int\limits_{\ext(V)}(l,\cdot)d\overline{\nu}=(l,y)
$$
for all $l\in\mathbb{R}^{N}$. Hence, $(\nu_{1},\nu_{2})$ belongs to~$U$
and
$$
\gathered
\nu_{1}(S_{1})+\nu_{2}(S_{2})=
\int\limits_{V_{1}}S_{1}d\overline{\nu}_{1}+
\int\limits_{V_{2}}S_{2}d\overline{\nu}_{2}
\\
=\int\limits_{V_{1}}\|\cdot\|_{S_{1}^{\circ}}d\overline{\nu}_{1}
+\int\limits_{V_{2}}\|\cdot\|_{S_{2}^{\circ}}d\overline{\nu}_{2}
=\overline{\nu}(\mathbb{1})=1=S_{1}\wedge S_{2}(y).
\endgathered
$$
By (5.1.1)
$$
\gathered
\sum\limits_{k=1}^{p}S_{1}\wedge S_{2}(x_{k})\le
\mu_{1}(S_{1})+\mu_{2}(S_{2})<\nu_{1}(S_{1})+\nu_{2}(S_{2})
\\
=S_{1}\wedge S_{2}(y)\le\sum\limits_{k=1}^{p}S_{1}\wedge S_{2}(x_{k}).
\endgathered
$$
We arrive at a contradiction, which means that
$(\mu_{1},\mu_{2})$ lies in~$\widetilde{U}$; i.~e. there are
measures $\nu_{1}$, $\nu_{2}$
such that
$\mu_{1}-\nu_{1}\in H^{*}$,
$\mu_{2}-\nu_{2}\in H^{*}$, and $(\nu_{1}$, $\nu_{2})\in U$.
Consider the representing points
$$
u_{1}:z\mapsto\nu_{1}(z);\quad
u_{2}:z\mapsto\nu_{2}(z)\quad (z\in\mathbbm{R}^{N}).
$$
Then $u_{1}+u_{2}=y$, and for $S\in H$
we have
$$
\mu_{1}(S)\ge\nu_{1}(S)\ge S(u_{1});\quad
\mu_{2}(S)\ge\nu_{2}(S)\ge S(u_{2}).
$$
Proceed by induction and apply the above process
to the measure $\mu_{2}$
and the nonzero point~$u_{2}$
(it is exactly the place where we invoke the simplification
of the beginning of the proof). We thus come to what was desired.
In case $u_{2}=0$, the sought decomposition
may be composed of the copies of the zero vectors.
The proof is complete.

By way of illustration of Theorem 5.1 we will
provide a~description for~$\pi(\Xi)$.

\subsection{}
\proclaim{Theorem.}
Let $H$ be a cone in~$\mathscr{V}_{N}$ and $H=\pi_\downarrow(H)$.
Assume that
$$
S_{y}:=
\bigwedge\limits_{S\in H;S\ne\{0\}}
\frac{S}{S(y)}
$$
is absorbing for every nonzero $y\in\mathbbm{R}^{N}$.
Then $\pi^\uparrow(H)$ is closed with respect to
$\wedge$. Moreover, and a nonzero $S$ in~$\mathscr{V}_{N}$ belongs to~$\pi^\uparrow(H)$
if and only if
$$
S=\bigwedge\limits_{x\ne 0}S(x)
\bigvee\limits_{S_{0}\in H}
\frac{S_{0}}{S_{0}(x)}
\leqno{(5.2.1)}
$$
\endproc

{\sc Proof.}
We have  already demonstrated that each $S\in\pi^\uparrow(H)$
may be written  as in~(5.2.1) (cp.~(4.2.2)).
Assume in turn that $S$ has the shape~(5.2.1).
By Theorem~4.1 we have to validate the implication
$$
\sum\limits_{k=1}^{n}S_{0}(x_{k})\ge S_{0}(y)\quad
\text{for all}\
S_{0}\in H\Longrightarrow\sum\limits_{k=1}^{n}S(x_{k})\ge S(y).
$$
Since $H=\pi^\uparrow(H)$, by Theorem~4.1 there are
vectors
$z_{1},\dots,z_{n}$ such that
$$
\gathered
\sum\limits_{k=1}^{n}z_{k}=y;
\\
S_{0}(x_{k})\ge S_{0}(z_{k})\quad (S_{0}\in H)
\endgathered
$$
Since  $S$ is represented as (5.2.1),
$S(x_{k})\ge S(z_{k})$. Hence,
$$
\sum\limits_{k=1}^{n}S(x_{k})\ge
\sum\limits_{k=1}^{n}S(z_{k})\ge
S\biggl(\,\sum\limits_{k=1}^{n}z_{k}\biggr)=S(y).
$$
Thus, $S\in^\pi\uparrow(H)$.

We are left with checking that $\pi^\uparrow(H)$
is closed under~$\wedge$.
By above, $S\in\pi^\uparrow(H)$
if and only if
$S(x)\ge S(y)$ for all $x,y\in \mathbbm{R}^{N}$ satisfying
$S_{0}(x)\ge S_{0}(y)$ for all $S_{0}\in  H$.

So, take $S_{1},S_{2}\in\pi^\uparrow(H)$
and assume that $S_{0}(x)\ge S_{0}(y)$ for all $S_{0}\in H$.

We are to compute $S_{1}\wedge S_{2}(x)$.
Arguing as in Theorem~5.1 and replacing the reference to
the Choquet Theorem to the Carath\'eodory Theorem,
find vectors $x_{1},\dots,x_{n}$
such that $\sum\limits_{k=1}^{n}x_{k}=x$ and
$$
S_{1}\wedge S_{2}(x)=\sum\limits_{k=1}^{t}S_{1}(x_{k})+
\sum\limits_{k=t+1}^{n}S_{2}(x_{k}).
$$
If $S_{0}\in H$ then
$$
\sum\limits_{k=1}^{n}S_{0}(x_{k})\ge
S_{0}\biggl(\,\sum\limits_{k=1}^{p}x_{k}\biggr)
=S_{0}(x)\ge S_{0}(y).
$$
Hence, by Theorem 5.1 there are vectors
$z_{1},\dots,z_{n}\in\mathbbm{R}^{N}$
such that $\sum\nolimits_{k=1}^{n}z_{k}=y$ and
$S_{0}(x_{k})\ge S_{0}(z_{k})$ for all $S_{0}\in H$
and $k:=1,\dots,n$. Thus, $S_{1}(x_{k})\ge S_{1}(z_{k})$ and
$S_{2}(x_{k})\ge S_{2}(z_{k})$.
Consequently,
$$
\gathered
S_{1}\wedge S_{2}(x)=\sum\limits_{k=1}^{t}S_{1}(x_{k})+
\sum\limits_{k=t+1}^{n}S_{2}(x_{k})\ge
\sum\limits_{k=1}^{t}S_{1}(z_{k})+
\sum\limits_{k=t+1}^{n}S_{2}(z_{k})
\\
\ge\sum\limits_{k=1}^{n}S_{1}\wedge S_{2}(z_{k})\ge
S_{1}\wedge S_{2}
\biggl(\,\sum\limits_{k=1}^{n}z_{k}\biggr)=S_{1}\wedge S_{2}(y).
\endgathered
$$
Therefore, $S_{1}\wedge S_{2}$ belongs to~$\pi^\uparrow(H)$,
which completes the proof.

\subsection{}
\proclaim{Corollary.}
Let $(S_{\xi})_{\xi\in\Xi}$ be a nondegenerate family
of balls. Then
$$
\pi(\Xi)=\pi^\uparrow(\pi_\downarrow(\Xi)).
$$
In this event a nonzero gauge~$S$
belongs to~$\pi(\Xi)$ if and only if
$$
S=\bigwedge\limits_{x\ne 0}S(x)\bigvee\limits_{S_{0}\in\pi_\downarrow(\Xi)}
\frac{S_{0}}{S_{0}(x)}.
$$
\endproc

{\sc Proof.}
Obviously, $\pi^\uparrow(\pi_\downarrow(\Xi))$ lies in~$\pi(\Xi)$.
Note now that
$$
S_{y}:=
\bigwedge\limits_{S_{0}\in\pi_\downarrow(\Xi);S_{0}\ne \{0\}}
\frac{S_{0}}{S_{0}(y)}\supset
\bigwedge\limits_{S_{0}\in {\rm M}(\Xi);S_{0}\ne \{0\}}
\frac{S_{0}}{S_{0}(y)}
$$
The family ${\rm M}(\Xi)$
is nondegenerate since so is
$(S_{\xi})_{\xi\in\Xi}$. Hence, $S_{y}$ is absorbing.
By Theorem 5.2 $\pi^\uparrow(\pi_\downarrow(\Xi))$
is closed under~$\wedge$, thus serving as
a~superset of~$\pi(\Xi)$.

\subsection{}
In study of the properties of gauges with are related to intersection,
we have actually used the accompanying representation
$$
\int\limits_{S_{N-1}}S_{1}\wedge S_{2}d\mu=
\inf\limits_{\mu_{1}+\mu_{2}
\underset{\mathbbm{R}^N}{\gg}\
\mu}
\biggl(\,\int\limits_{S_{N-1}}S_{1}d\mu_{1}+
\int\limits_{S_{N-1}}S_{2}d\mu_{2}\biggr),
\leqno{(5.4.1)}
$$
which generalizes the standard formula for
the {\it infimal convolution\/} $\square$, a~routine operation
of convex analysis:
$$
S_{1}\wedge S_{2}=S_{1}\square\, S_{2}.
$$
It is an easy matter to see the lattice-theoretic
provenance of~(5.4.1). Some slightly annoying subtlety of
the general case which was obviated by finite dimensionality
is connected with the fact the infimum of
abstract convex elements in the lattice of these elements
is just a~partial superlinear operator.

{\bf Acknowledgement.} The main results of this article stem
from our joint work with Alex Rubinov by the mid 1970s.
I gratefully emphasize his creative contribution to all areas
of abstract convexity  we  had been exploring those happy years.

\bibliographystyle{plain}

\enddocument